\documentclass[11pt,final]{amsart}
\usepackage{amsmath}
\usepackage{amssymb}

\theoremstyle{plain}
\newtheorem*{theorem}{Theorem}

\newtheorem{corollary}{Corollary}

\theoremstyle{remark}
\newtheorem{remark}{Remark}

\newcommand\C{\mathbb C}
\renewcommand\ker{\textrm{ker }}
\newcommand\lam{\lambda}
\newcommand\la{\langle}
\newcommand\ot{\otimes}

\newcommand\ra{\rangle}
\newcommand\rng{\textrm{rng }}
\newcommand\tr{\textrm{tr }}
\newcommand\ul{\underline}

\begin{document}
\title[]{Generalization of Wigner's Unitary-Antiunitary Theorem for
Indefinite Inner Product Spaces}
\author{Lajos Moln\'ar}
\dedicatory{Dedicated to M\'ad}
\address{Institute of Mathematics and Informatics,
         Lajos Kossuth University,
         4010 Debrecen, P.O.Box 12, Hungary}
\email{molnarl@math.klte.hu}
\maketitle
\begin{abstract}
We present a generalization of Wigner's unitary-antiunitary theorem
for pairs of ray transformations. As a particular case, we get
a new Wigner-type theorem for non-Hermitian indefinite inner product
spaces.
\end{abstract}
\maketitle

\vskip 24pt
\noindent
The classical Wigner unitary-antiunitary theorem plays a fundamental
role in the foundations of quantum mechanics and it also has deep
connections with the theory of projective spaces. It states that every
ray transformation (see below) on a Hilbert space which preserves
the transition probabilities can be lifted to a (linear) unitary or
a (conjugate-linear) antiunitary operator on $H$ (see
\cite{Bar,Rat,ShA1}).
So, Wigner's result concerns definite inner product spaces.
On the other hand, it has become quite clear by now that the indefinite
inner product spaces might be even more useful for the discussion of
several physical problems. For example, this is the case in relation to
the divergence problem in quantum field theory, or when one wants to
preserve some basic properties of the field like relativistic
covariance and locality (see the introduction of \cite{BMS}).
This raises the need to study Wigner's theorem in the "indefinite"
setting as well. Previous results in this direction were presented in
\cite{Bau,BMS,Bro}.
The aim of this paper is to contribute to this study by giving a
very general Wigner-type theorem which involves not one but two
ray transformations and then apply it to get a generalization of
Wigner's theorem for indefinite inner product spaces. The main
difference which distiguishes our result from the previous ones is
that we do not assume even that the indefinite inner product
under consideration
is Hermitian. What allows us to reach this result
is that we refine our algebraic approach to Wigner's theorem
which has already been proved to be fruitful in our recent papers
\cite{MolJAMS,MolJMP}.
The main feature of this approach is that instead of manipulating in the
underlying space, we push the problem to an operator algebra over
our space and apply some classical results from pure ring theory.
Hence,
our method is completely different from those used previously in the
papers dealing with Wigner's theorem in indefinite inner product spaces.

Let us fix the definitions and notation that we shall use throughout.
In what follows, let $H$ be a Hilbert space.
Given a vector $x\in H$, the set of all vectors of the form $\lam x$
with
$\lam \in \C$, $|\lam |=1$ is called the ray associated to $x$ and it
is denoted by $\ul x$. For any $x, y\in H$ we define
\[
\ul x \cdot \ul y=|\la x,y\ra |.
\]
The notation $\ul H$ stands for the set of all rays in $H$.
The algebra of all bounded linear operators on $H$ is denoted by $B(H)$,
and $F(H)$ stands for the ideal of all finite rank operators in $B(H)$.
If $x,y\in H$ are arbitrary vectors, then $x\ot y$ is an element of
$F(H)$ which is defined by $(x\ot y)z=\la z,y\ra x$ $(z\in H)$.
A linear map $\phi:\mathcal A \to \mathcal B$ between the algebras
$\mathcal A$ and $\mathcal B$ is called a Jordan homomorphism if
\[
\phi(x^2)=\phi(x)^2 \qquad (x\in \mathcal A).
\]

Our main result which follows presents a Wigner-type result for
pairs of ray transformations.

\begin{theorem}
Let $H$ be a complex Hilbert space of dimension at least 3.
Let $T, S: \ul H \to \ul H$ be bijective transformations with the
property that
\[
T\ul x \cdot S\ul y = \ul x \cdot \ul y \qquad (\ul x, \ul y\in \ul H).
\]
Then there are bounded invertible either both linear or both
conjugate-linear operators $U,V:H \to H$ such that $V={U^*}^{-1}$ and
\[
T\ul x ={\underline{Ux}}, \qquad
S\ul x ={\underline{Vx}} \qquad (x\in H).
\]
\end{theorem}

\begin{proof}
For every $x\in H$ pick a vector from $T \ul x$. In that way we get a
function, which will be denoted by the same symbol $T$, from $H$ into
itself with the property that for every vector $y\in H$, there exists a
vector $x\in H$ such that $y=\lam Tx$ for some $\lam \in \C$ of modulus
1. Let us do the same with the other transformation $S$. Clearly, we
have
\[
|\la Tx, Sy\ra |=|\la x,y\ra | \qquad (x,y \in H).
\]
Obviously,
for every unit vector $x\in H$ we can choose a scalar $\lam_x $ with
$|\lam_x |=1$ such that $\lam_x \la Tx, Sx\ra =1$.
By the properties of our original transformation $T$, we can clearly
suppose that here in fact we have $\la Tx,Sx\ra =1$.
We define a function $\mu$ on the set $P_f(H)$ of all finite rank
projections (self-adjoint idempotents) on $H$ as follows. If $P\in
P_f(H)$, then there are pairwise
orthogonal unit vectors $x_1, \ldots ,x_n\in H$ such that $P=x_1\ot x_1
+\ldots +x_n\ot x_n$. We set
\[
\mu(P)= Tx_1 \ot Sx_1 +\ldots + Tx_n \ot Sx_n.
\]
Apparently, the operators
$Tx_1 \ot Sx_1,\ldots , Tx_n \ot Sx_n$ are
pairwise orthogonal rank-one idempotents (two idempotents $P,Q$ are
said to be orthogonal if $PQ=QP=0$). Hence, $\mu(P)$ is a rank-$n$
idempotent. We have to check that $\mu$ is well-defined. This follows
from the following observation. We have
\[
\rng (\sum_{k=1}^n  Tx_k \ot Sx_k)=[Tx_1, \ldots, Tx_n]
\]
and
\[
\ker (\sum_{k=1}^n Tx_k \ot Sx_k)=[Sx_1, \ldots, Sx_n]^\perp,
\]
where $[.]$ denotes generated subspace. Now, suppose that the pairwise
orthogonal unit vectors $x_1', \ldots, x_n'$ generate the same subspace
as
$x_1,\ldots, x_n$ do. Let $y \in H$. Then there exist a vector $x\in H$
and a scalar $\lam $ of modulus 1 such that $y=\lam Sx$. We have
\[
y\perp [Tx_1, \ldots, Tx_n] \Leftrightarrow
Sx\perp [Tx_1, \ldots, Tx_n] \Leftrightarrow
\]
\begin{equation}\label{E:star1}
x\perp [x_1, \ldots, x_n] \Leftrightarrow
x\perp [x_1', \ldots, x_n'] \Leftrightarrow
\end{equation}
\[
Sx\perp [Tx_1', \ldots, Tx_n'] \Leftrightarrow
y\perp [Tx_1', \ldots, Tx_n'] .
\]
This shows that the range of
$\sum_{k=1}^n  Tx_k \ot Sx_k$
is the same as that of
$\sum_{k=1}^n  Tx_k' \ot Sx_k'$. The same applies for the
kernels. Since the idempotents are determined by their ranges and
kernels, this proves that $\mu$ is well-defined.
It is now clear that $\mu$ is an orthoadditive measure on $P_f(H)$.
We show that $\mu$ is bounded on the set $P_1(H)$ of all rank-one
projections which is equivalent to
\[
\sup_{\|x\|=1} \| Tx\| \| Sx\| < \infty.
\]
Suppose, on the contrary, that there is a sequence $(u_n)$ of unit
vectors in $H$ for which $\| Tu_n\| \| Su_n\| \longrightarrow \infty$.
Since $(u_n)$ is bounded, it has a subsequence
$(u_{k_n})$ weakly converging to a vector, say, $u\in H$. We have
\[
|\la Tu_{k_n} ,Sv\ra| =
|\la u_{k_n} ,v\ra| \longrightarrow |\la u,v\ra |.
\]
Since this holds for every $v\in H$, we deduce that $(Tu_{k_n})$ is
weakly bounded which implies that it is in fact norm-bounded. The same
argument applies in relation to $S$. Hence, we obtain that $(u_n)$ has a
subsequence $(u_{l_n})$ such that $\| Tu_{l_n}\|, \| Su_{l_n}\|$ are
bounded which is a contradiction. Consequently, $\mu$ is
bounded on $P_1(H)$.

By Gleason's theorem $\mu$ can be extended to a Jordan homomorphism of
$F(H)$. In fact, if $A\in F(H)$ is self-adjoint, then there are
finite-rank projections $P_1, \ldots, P_n$ (here, we do not require
that they are pairwise orthogonal) and scalars $\lam_1, \ldots, \lam_n$
such that $A=\lam_1 P_1+\ldots +\lam_n P_n$. Let
\[
\phi(A) =\lam_1 \mu(P_1) +\ldots +\lam_n \mu(P_n).
\]
Consider a finite dimensional subspace $H_0$ of $H$ with dimension at
least 3 which contains all the subspaces $\rng A, \ker A^\perp, \rng
P_1, \ldots, \rng P_n$.
Since $\mu$ is bounded on $P_1(H_0)$, by
the variation \cite[Theorem 3.2.16]{Dvur} of Gleason's theorem, for
every $x,y\in H$ there is an operator $T_{xy}$ on $H_0$ such that
\[
\la \lam_1 \mu(P_1) +\ldots +\lam_n\mu(P_n)x, y\ra=
\lam_1 \la \mu(P_1)x,y\ra  +\ldots +\lam_n \la \mu(P_n)x, y\ra=
\]
\[
\lam_1 \tr (P_1 T_{xy})  +\ldots +\lam_n \tr (P_n T_{xy})=
\tr (A T_{xy}).
\]
We now easily obtain that $\phi$ is well-defined and real-linear on the
set of all self-adjoint
finite rank operators. If $A\in F(H)$ is arbitrary, then there exist
self-adjoint finite rank operators $A_1, A_2$ such that $A=A_1+iA_2$.
Define $\phi(A)=\phi(A_1)+i\phi(A_2)$. Clearly, $\phi$ is a linear map
on $F(H)$ which sends projections to idempotents. It is a standard
algebraic argument to verify that $\phi$ is then a Jordan homomorphism
(see, for example, the proof of \cite[Theorem 2]{MolStud1}).
Since $F(H)$ is a locally matrix ring, we can apply a classical theorem
of Jacobson and Rickart. By \cite[Theorem 8]{JR} we obtain that $\phi$
can
be written as $\phi=\phi_1+\phi_2$, where $\phi_1$ is a homomorphism and
$\phi_2$ is an antihomomorphism. Since $\phi(P)$ is a rank-one
idempotent and $\phi_1(P)$, $\phi_2(P)$ are
idempotents, we infer from $\phi(P)=\phi_1(P)+\phi_2(P)$ that either
$\phi_1(P)=0$
or $\phi_2(P)=0$. Since the ring $F(H)$ is simple, we obtain that either
$\phi_1=0$ or $\phi_2=0$. Therefore, $\phi$ is either a homomorphism or
an antihomomorphism.
Without loss of generality we can assume that $\phi$ is a homomorphism.
We assert that $\phi$ is rank-preserving. Let $A\in F(H)$ be a
rank-$n$ operator. Then there is a rank-$n$ projection $P$ such that
$PA=A$. The rank of $\phi(P)$ is also $n$.
We have $ \phi(A)= \phi(PA)=\phi(P) \phi(A)$
which proves that $ \phi(A)$ is of rank at most $n$. If $Q$ is
any rank-$n$ projection, then there are finite rank operators $U,V$
such that $Q=UAV$. Since $\phi(Q)= \phi(U) \phi (A)
\phi(V)$ and the rank of $\phi(Q)$ is $n$, it follows that the rank of
$ \phi(A)$ is at least $n$. Therefore, $ \phi$ is rank-preserving. We
now refer to Hou's work \cite{Hou} on the form of linear
rank preservers on operator algebras.
It follows from the argument leading to \cite[Theorem 1.2]{Hou}
that there are linear
operators $U,V$ on $H$ such that $ \phi$ is of the form
\begin{equation}\label{E:2gleasegy}
\phi(x\ot y)=(Ux)\ot (Vy) \qquad (x,y \in H)
\end{equation}
(recall that we have assumed that $ \phi$ is a homomorphism).
If $x\in H$ is a unit vector, then we have
$Tx\ot Sx=\phi(x\ot x)=Ux\ot Vx$. Taking traces, we obtain
$1=\la Tx,Sx\ra=\la Ux,Vx\ra$. Since this holds for every unit vector
$x$, by the linearity of $U,V$, using polarization we get that
\begin{equation}\label{E:star2}
\la Ux,Vy\ra=\la x,y\ra \qquad (x,y\in H).
\end{equation}
We assert that $U,V$ are surjective. Consider, for example, the case of
$U$.
Let $0\neq x\in H$ be any vector and let $0\neq \lambda\in \C$ be any
scalar. It is easy to see that
$[Tx]^\perp =[T(\lambda x)]^\perp$ (see (\ref{E:star1})).
Therefore, $T(\lam x)= \lambda' Tx $ with some scalar $\lam'$.
Denote $x_e=x/\| x\|$. We compute
\[
Ux\ot Vx=\|x\|^2 Ux_e \ot Vx_e=
\|x\|^2 \phi(x_e\ot x_e)=
\| x\|^2 Tx_e\ot Sx_e.
\]
This gives us that $Tx_e \in [Ux]$. But $Tx$ is in the one-dimensional
subspace generated by $Tx_e$. So, we have
\begin{equation}\label{E:star3}
Tx\in [Ux].
\end{equation}
Since $\rng U$ is a linear subspace of $H$ and $T$ is "almost"
surjective, we obtain the surjectivity of $U$.
Similar argument applies to $V$.
We next show that $U,V$ are bounded.
Let $(x_n)$ be a sequence converging to 0 and let $y\in H$ be
such that $Ux_n \to y$. If $x\in H$ is arbitrary, then we have
\[
\la Ux_n, Vx\ra =\la x_n ,x\ra \longrightarrow 0.
\]
Since $V$ is surjective, we obtain that $(Ux_n)$ weakly
converges to 0. It follows that $y=0$. By the closed graph theorem we
deduce that $U$ is bounded. Similar argument proves the boundedness
of $V$. It follows from (\ref{E:star2}) that $V^*U=I$. This gives us
that $U$ is injective. Therefore, $U$ and $V$ are invertible and
$V={U^*}^{-1}$.

By (\ref{E:star3}) and the similar relation $Sx\in [Vx]$ $(x\in H)$,
there are functions $\varphi,\psi :H \to \C$ such that
\[
Tx =\varphi(x)Ux, \qquad Sx=\psi(x) Vx \qquad (x\in H).
\]
We have
\[
|\varphi(x)||\psi(y)||\la x,y\ra|=
|\varphi(x)||\psi(y)||\la Ux,Vy\ra|=
|\la Tx,Sy\ra|=
|\la x,y\ra|,
\]
that is, $|\varphi(x)||\psi(y)|=1$ if $\la x,y\ra\neq 0$. This easily
implies that $|\varphi|$ and $|\psi|$ are both constant. Multiplying
$U,V$, $\varphi,\psi$ by suitable constants, we obtain the statement of
the theorem. The proof is complete.
\end{proof}

In the following corollary of our theorem we give a generalization of
Wigner's theorem for the indefinite
inner product space generated by any invertible operator $A\in B(H)$.
Since we do not assume that $A$ is self-adjoint,
this result can, in some sense, be considered as a generalization of
the results in \cite{Bau,BMS}.

\begin{corollary}\label{C:wigelso}
Let $H$ be a complex Hilbert space with $\dim H\geq 3$ and let $A\in
B(H)$ be invertible.
For any $x,y\in H$ define $\ul x \cdot_A \ul y=|\la Ax, y\ra|$.
Let $T:\ul H \to \ul H$ be a bijective transformation such that
\[
T\ul x\cdot_A Ty=\ul x\cdot_A \ul y \qquad (x,y\in H).
\]
Then there is a bounded invertible either linear or conjugate-linear
operator $U$ on $H$ with $U^*AU=\epsilon A$ for some scalar
$\epsilon$ of modulus 1 such that
\[
T\ul x ={\underline{Ux}} \qquad (x\in H).
\]
\end{corollary}

\begin{proof}
Just as in the proof of our theorem above, we can define an "almost"
surjective map (that is, which has values in every ray) on the
underlying Hilbert space $H$ denoted by the same symbol $T$ such that
\[
|\la ATx,Ty\ra|=|\la Ax,y\ra | \qquad (x,y\in H).
\]
Set $S=ATA^{-1}$. The proof of our theorem now applies and we find that
there is a bounded invertible either linear or conjugate-linear operator
$U$ on $H$ and a scalar function $\varphi: H \to \C$ such that
$Tx =\varphi(x) Ux$ $(x \in H)$. Since
\begin{equation}\label{E:sutty}
|\varphi(x)||\varphi(y)| |\la AUx, Uy\ra |=
|\la ATx,Ty\ra |=
|\la Ax, y\ra | \qquad (x,y\in H),
\end{equation}
it follows that $[U^*AUx]^\perp= [Ax]^\perp$ for every $x\in H$.
Therefore, the linear operators $U^*AU$ and $A$
are locally linearly dependent which means that $U^*AUx$ and $Ax$ are
linearly dependent for every $x\in H$. Since none of the operators
$U^*AU$ and $A$ is of rank 1,
by \cite[Lemma 3]{Hou2} we obtain that there is a scalar $c$ such that
$U^*AU=cA$. Let $x,y\in H$ be arbitrary nonzero vectors. Pick $z\in H$
such that $\la Ax,z\ra, \la Ay,z\ra \neq 0$.
From (\ref{E:sutty}) we now infer that
\[
|\varphi(x)||\varphi(z)| |c|=1, \qquad
|\varphi(y)||\varphi(z)| |c|=1.
\]
This shows that $|\varphi|$ is constant. If $d$ denotes this constant,
then we have $d^2|c|=1$. Let $\epsilon=d^2 c$. Then $\epsilon$ is of
modulus 1 and we have
\[
(d\epsilon U)^*A (d\epsilon U)=d^2 U^*AU=d^2c A=\epsilon A.
\]
Consider the factorization
\[
Tx=\biggl(\frac{1}{d\epsilon} \varphi(x)\biggr) (d\epsilon U).
\]
Since $\frac{1}{d\epsilon} \varphi(x)$ is of modulus 1, the proof is
complete.
\end{proof}

In the finite dimensional case, Corollary~\ref{C:wigelso} can be
reformulated in the following way.

\begin{corollary}
Let $H$ be a finite dimensional complex Hilbert space with $\dim
H\geq 3$. Let $B:H\times H \to \C$ be a sesquiliner form which
is non-degenerate in the sense that $B(x,y)=0$ $(y\in H)$ implies $x=0$.
Define $\ul x\cdot_B \ul y= |B(x,y)|$ $(x,y \in H)$.
Let $T:\ul H \to \ul H$ be a bijective transformation such that
\[
T\ul x\cdot_B Ty=\ul x\cdot_B \ul y \qquad (x,y\in H).
\]
Then either there is an invertible linear
operator $U$ on $H$ such that $B(Ux,Uy)=\epsilon B(x,y)$ $(x,y\in H)$
for some scalar $\epsilon$ of modulus 1 and
\[
T\ul x ={\underline{Ux}} \qquad (x\in H),
\]
or there is an invertible conjugate-linear
operator $U'$ on $H$ such that $\overline{B(U'x,U'y)}=\epsilon' B(x,y)$
$(x,y\in H)$ for some scalar $\epsilon'$ of modulus 1 and
\[
T\ul x ={\underline{U'x}} \qquad (x\in H).
\]
\end{corollary}

\begin{proof}
Since $H$ is finite dimensional, it is easy to see that there exists
an invertible linear operator $A$ on $H$ such that
$B(x,y)=\la Ax,y\ra$ $(x,y\in H)$. Now, Corollary~\ref{C:wigelso}
applies.
\end{proof}

\begin{remark}
Our results are valid in real Hilbert spaces as well.
In order to see it, we must refine the argument we
have presented in the complex case. Namely, one can follow the
argument
that has been applied in the proof of \cite[Theorem 3]{MolJAMS}.
Observe that in the papers \cite{Bau,BMS}
the authors considered only complex spaces.
\end{remark}

\begin{center}
\bf{Acknowledgements}
\end{center}
          This research was supported from the following sources:
          (1) Hungarian National Foundation for Scientific Research
          (OTKA), Grant No. T--030082 F--019322,
          (2) A grant from the Ministry of Education, Hungary, Reg.
          No. FKFP 0304/1997.


\end{document}